\documentclass{amsart}
\usepackage{amsmath,amssymb}
\usepackage{yhmath}

\newcommand{\bdis}{\begin{displaymath}}
\newcommand{\edis}{\end{displaymath}}
\newcommand{\be}{\begin{equation}}
\newcommand{\ee}{\end{equation}}
\newcommand{\mbb}{\mathbb}
\newcommand{\mcal}{\mathcal}

\newcommand{\vp}{\varphi}

\newcommand{\zf}{\zeta\left(\frac{1}{2}+it\right)}
 
\newcommand{\zfn}{\zeta\left(\frac{1}{2}+it_\nu\right)}   
\newcommand{\zfnn}{\zeta\left(\frac{1}{2}+it_{\nu+1}\right)}

\DeclareMathOperator{\re}{Re}

\theoremstyle{definition}

\theoremstyle{remark}
\newtheorem{remark}[]{Remark}

\newtheorem*{mydef11}{{\bf Theorem 1}}

\newtheorem*{mydef12}{{\bf Theorem 2}}

\newtheorem*{mydef41}{{\bf Corollary 1}}

\newtheorem*{mydef51}{{\bf Lemma 1}}

\newtheorem*{mydef52}{{\bf Lemma 2}}

\newtheorem*{mydef53}{{\bf Lemma 3}}

\newtheorem*{mydef81}{{\bf Property 1}}

\newtheorem*{mydef91}{{\bf Formula 1}}

\newtheorem*{mydef92}{{\bf Formula 2}}

\newtheorem*{mydef93}{{\bf Formula 3}}

\numberwithin{equation}{section}

\begin{document}

\title[Jacob's ladders, almost exact decomposition \dots]{Jacob's ladders, almost exact decomposition of certain increments of the Hardy-Littlewood integral (1918) by means of the Raabe's integral and the thirteenth equivalent of the Fermat-Wiles theorem}

\author{Jan Moser}

\address{Department of Mathematical Analysis and Numerical Mathematics, Comenius University, Mlynska Dolina M105, 842 48 Bratislava, SLOVAKIA}

\email{jan.mozer@fmph.uniba.sk}


\begin{abstract}
In this paper we use our theory of Jacob's ladders on the Raabe's integral to obtain: (i) The thirteenth equivalent of the Fermat-Wiles theorem, as well as (ii) almost exact decomposition of certain elements of continuum set of increments of the Hardy-Littlewood integral. 
\end{abstract}
\maketitle

\section{Introduction} 

\subsection{} 

Let us remind that we have introduced\footnote{See \cite{8}, (4.6).} a new kind of functional, namely 
\be \label{1.1} 
\lim_{\tau\to\infty}\frac{1}{\tau}\int_{\frac{x}{1-c}\tau}^{[\frac{x}{1-c}\tau]^1}\left|\zf\right|^2{\rm d}t=x
\ee 
for every fixed $x>0$, where 
\be \label{1.2} 
\left[\frac{x}{1-c}\tau\right]^1=\vp_1^{-1}\left(\frac{x}{1-c}\tau\right). 
\ee 
Every fixed ray 
\be \label{1.3} 
y(\tau;x)=\frac{x}{1-c}\tau,\ \tau\to\infty 
\ee  
is associated with given positive $x$ defined by the slope 
\be \label{1.4} 
0<\arctan\frac{x}{1-c}<\frac{\pi}{2}. 
\ee  

\subsection{} 

Our functional (\ref{1.1}) represents the basis for the construction of completely new set of $\zeta$-equivalents of the Fermat-Wiles theorem, see the paper \cite{8}. 

We have also used the functional (\ref{1.1}) for construction of another set of equivalents of the Fermat-Wiles theorem. These are expressed by means of the following functions and corresponding formulas\footnote{All these constructions can be found in our papers \cite{9} -- \cite{12}.}: 
\begin{itemize}
	\item Dirichlet's function 
	\be \label{1.5} 
	D(x)=\sum_{n\leq x}d(n), 
	\ee  
	\item Titchmarsh's functions
	\be \label{1.6} 
	\mcal{T}_1(x)=\sum_{t_\nu\leq x}\zfn, 
	\ee 
	\be \label{1.7} 
	\mcal{T}_2(x)=\sum_{t_\nu\leq x}\zfn\zfnn, 
	\ee 
	\item Stirling's formula for the Euler's $\Gamma$-function 
	\be \label{1.8} 
	\Gamma(x)=x^{x-\frac 12}e^{-x}\sqrt{2\pi}e^{\frac{\Theta}{12x}},\ 0<\Theta<1, 
	\ee 
	\item logarithmic modification of the Hardy-Littlewood integral 
	\be \label{1.9} 
	J_1(T)=\int_0^T\left\{\ln t-\left|\zf\right|^2\right\}{\rm d}t, 
	\ee  
	\item Selberg's formula (1946) for the function 
	\be \label{1.10} 
	S_1(t)=\frac{1}{\pi}\int_0^t\arg\zf{\rm d}t
	\ee 
	in a linear combination with (\ref{1.1}). 
\end{itemize} 

\subsection{} 

In this paper we use our theory in the case the Raabe's intergral 
\be \label{1.11} 
\int_a^{a+1}\ln\Gamma(t){\rm d}t=a\ln a - a+\ln\sqrt{2\pi},\ a>0 
\ee 
and we obtain the following result. The condition 
\be \label{1.12} 
\begin{split}
& \lim_{\tau\to\infty}\frac{1}{\tau}\left\{
\int_{\frac{x^n+y^n}{z^n}\frac{\tau}{1-c}}^{\frac{x^n+y^n}{z^n}\frac{\tau}{1-c}+1}\ln\Gamma(t){\rm d}t-\int_{\vp_1(\frac{x^n+y^n}{z^n}\frac{\tau}{1-c})}^{\vp_1(\frac{x^n+y^n}{z^n}\frac{\tau}{1-c})+1}\ln\Gamma(t){\rm d}t
\right\}\not=1
\end{split}
\ee 
on the class of all Fermat's rationals 
\be \label{1.13} 
\frac{x^n+y^n}{z^n},\ x,y,z,n\in\mbb{N},\ n\geq 3
\ee 
represents the thirteenth equivalent of the Fermat-Wiles theorem. 

\begin{remark}
Of course, it is sufficient to localize the condition (\ref{1.12}) in the way 
\bdis 
1-\epsilon<\frac{x^n+y^n}{z^n}<1+\epsilon 
\edis 
for every small, positive and fixed $\epsilon$, comp. \cite{7}, Sec. 7. 
\end{remark} 

\subsection{} 

Next, we obtain the following result: For the increments 
\bdis 
\begin{split}
& \int_{\overset{r}{T}}^{\overset{r+1}{T}}\left|\zf\right|^2{\rm d}t,\ \overset{r}{T}=\overset{r}{T}(T), \\ 
& r=1,\dots,k,\ k\in\mbb{N}
\end{split}
\edis 
there are the following almost exact decompositions 
\be \label{1.14} 
\begin{split}
& \int_{\overset{r}{T}}^{\overset{r+1}{T}}\left|\zf\right|^2{\rm d}t = \\ 
& \int_{\overset{r}{T}}^{\overset{r}{T}+1}\ln\Gamma(t){\rm d}t-\int_{\overset{r-1}{T}}^{\overset{r-1}{T}+1}\ln\Gamma(t){\rm d}t- \\ 
& (\ln 2\pi-1-c)(\overset{r}{T}-\overset{r-1}{T})+\mcal{O}\left(\frac{\ln T}{T}\right),\ T\to\infty 
\end{split}
\ee 
for every fixed $k$. 

\begin{remark}
It is clear that the almost exact decompositions (\ref{1.14}) are completely inaccessible by the Hardy-Littlewood-Ingham formula\footnote{Comp. \cite{9}, (3.4) and (3.5).}
\be \label{1.15} 
\int_0^T\left|\zf\right|^2{\rm d}t=T\ln T-(1+\ln 2\pi-2c)T+R(T)
\ee 
with the error term 
\be \label{1.16} 
R(T)=\mcal{O}(T^{a+\delta}),\ \frac 14\leq a\leq \frac 13, 
\ee 
where the value $\frac 14$ is fixed by the fundamental Good's $\Omega$-theorem (1977) and $\delta>0$ is a small number. 
\end{remark}

\subsection{} 

In this paper we use the following notions of our works \cite{2} -- \cite{5}: 
\begin{itemize}
	\item[{\tt (a)}] Jacob's ladder $\vp_1(T)$, 
	\item[{\tt (b)}] direct iterations of Jacob's ladders 
	\bdis 
	\begin{split}
		& \vp_1^0(t)=t,\ \vp_1^1(t)=\vp_1(t),\ \vp_1^2(t)=\vp_1(\vp_1(t)),\dots , \\ 
		& \vp_1^k(t)=\vp_1(\vp_1^{k-1}(t))
	\end{split}
	\edis 
	for every fixed natural number $k$, 
	\item[{\tt (c)}] reverse iterations of Jacob's ladders 
	\be \label{1.17}  
	\begin{split}
		& \vp_1^{-1}(T)=\overset{1}{T},\ \vp_1^{-2}(T)=\vp_1^{-1}(\overset{1}{T})=\overset{2}{T},\dots, \\ 
		& \vp_1^{-r}(T)=\vp_1^{-1}(\overset{r-1}{T})=\overset{r}{T},\ r=1,\dots,k, 
	\end{split} 
	\ee   
	where, for example, 
	\be \label{1.18} 
	\vp_1(\overset{r}{T})=\overset{r-1}{T}
	\ee  
	for every fixed $k\in\mbb{N}$ and every sufficiently big $T>0$. We also use the properties of the reverse iterations listed below.  
	\be \label{1.19}
	\begin{split} 
		& \overset{r}{T}-\overset{r-1}{T}\sim(1-c)\pi(\overset{r}{T});\ \pi(\overset{r}{T})\sim\frac{\overset{r}{T}}{\ln \overset{r}{T}},\ r=1,\dots,k,\ T\to\infty, \\ 
		& \overset{0}{T}=T<\overset{1}{T}(T)<\overset{2}{T}(T)<\dots<\overset{k}{T}(T), \\ 
		& T\sim \overset{1}{T}\sim \overset{2}{T}\sim \dots\sim \overset{k}{T},\ T\to\infty.   
	\end{split}
	\ee  
\end{itemize} 

\begin{remark}
	The asymptotic behaviour of the points 
	\bdis 
	\{T,\overset{1}{T},\dots,\overset{k}{T}\}
	\edis  
	is as follows: at $T\to\infty$ these points recede unboundedly each from other and all together are receding to infinity. Hence, the set of these points behaves at $T\to\infty$ as one-dimensional Friedmann-Hubble expanding Universe. 
\end{remark}  

\section{Riemann's zeta-function and \emph{proliferation} of every $L_2$-orthogonal system} 

Let us remind that the construction of the set of thirteenth equivalents of the Fermat-Wiles theorem mentioned in the subsections 1.2 and 1.3 gives simultaneously the set of the points of contact between completely different mathematical objects. They are namely: 
\begin{itemize}
	\item[(A)] Our functional (\ref{1.1}) and the set of formulas (\ref{1.5}) -- (\ref{1.11}); 
	\item[(B)] The Fermat-Wiles theorem and the set of formulae mentioned in (A). 
\end{itemize} 

For completeness, we remind also our result concerning \emph{the proliferation}, see \cite{6}, that represents the point of contact between the Riemann's function 
\bdis 
\zf 
\edis  
and the theory of $L_2$-orthogonal systems. We shortly recapitulate that result in the text below.  

\subsection{} 

We have introduced the generating vector operator $\hat{G}$ acting on the class of all $L_2$-orthogonal systems 
\be \label{2.1} 
\{f_n(t)\}_{n=0}^\infty,\ t\in [a,a+2l],\ a\in\mbb{R},\ l>0 
\ee 
as 
\be \label{2.2} 
\begin{split}
& \{f_n(t)\}_{n=0}^\infty\xrightarrow{\hat{G}}\{f_n^{p1}(t)\}_{n=0}^\infty\xrightarrow{\hat{G}}\{f_n^{p_1,p_2}(t)\}_{n=0}^\infty\xrightarrow{\hat{G}} \dots \\ 
& \{f_n^{p_1,p_2,\dots,p_s}(t)\}_{n=0}^\infty,\ p_1,\dots,p_s=1,\dots k 
\end{split}
\ee 
for every fixed $k,s\in\mbb{N}$ with explicit formulae\footnote{See \cite{6}, (2.19).} for 
\bdis 
f_n^{p_1,p_2,\dots,p_s}(t). 
\edis  

\subsection{} 

In the case of Legendre's orthogonal system 
\be \label{2.3} 
\{ P_n(t)\}_{n=0}^\infty,\ t\in [-1,1] 
\ee 
the operator $\hat{G}$ produces, for example, the third generation as follows 
\be \label{2.4} 
\begin{split}
& P_n^{p_1,p_2,p_3}(t)=P_n(u_{p_1}(u_{p_2}(u_{p_3}(t))))\times\prod_{r=0}^{p_1-1}|\tilde{Z}(v_{p_1}^r(u_{p_2}(u_{p_3}(t))))|\times \\ 
& \prod_{r=0}^{p_2-1}|\tilde{Z}(v_{p_2}^r(u_{p_3}(t)))|\times\prod_{r=0}^{p_3-1}|\tilde{Z}(v_{p_3}^r(t))|, \\ 
& p_1,p_2,p_3=1,\dots,k,\ t\in [-1,1],\ a=-1,\ l=1, 
\end{split}
\ee 
where 
\be \label{2.5} 
u_{p_i}(t)=\vp_1^{p_i}\left(\frac{\overset{p_i}{\wideparen{T+2}}-\overset{p_i}{T}}{2}(t+1)+\overset{p_i}{T}\right)-T-1,\ i=1,2,3 
\ee  
are the automorphisms on $[-1,1]$, and 
\be \label{2.6} 
\begin{split}
& v_{p_i}^r(t)=\vp_1^{r}\left(\frac{\overset{p_i}{\wideparen{T+2}}-\overset{p_i}{T}}{2}(t+1)+\overset{p_i}{T}\right),\ r=0,1,\dots,p_i-1, \\ 
& t\in [-1,1] \ \Rightarrow \ u_{p_i}(t)\in [-1,1] \ \wedge \ v_{p_i}^r(t)\in [\overset{p_i-r}{T},\overset{p_i-r}{\wideparen{T+2}}]. 
\end{split}
\ee  

\begin{mydef81}
\begin{itemize}
	\item[(a)] Every member of every $L_2$-orthogonal system 
	\be \label{2.7} 
	\{P_n^{p_1,p_2,p_3}(t)\}_{n=0}^\infty,\ t\in [-1,1],\ p_1,p_2,p_3=1,\dots, k
	\ee 
	contains the function 
	\bdis 
	\left|\zf\right|_{t=\tau}
	\edis 
	for corresponding $\tau$ since\footnote{See \cite{3}, (9.1), (9.2).} 
	\be \label{2.8} 
	|\tilde{Z}(t)|=\sqrt{\frac{{\rm d}\vp_1(t)}{{\rm d}t}}=\frac{1+o(1)}{\sqrt{\ln t}}\left|\zf\right|,\ t\to\infty. 
	\ee 
	\item[(b)] The Property (a) holds true due to Theorem of the paper \cite{6} for every generation 
	\bdis 
	\{f_n^{p_1,p_2,\dots,p_s}(t)\}_{n=0}^\infty,\ t\in [a,a+2l],\ s\in\mbb{N}. 
	\edis 
\end{itemize}
\end{mydef81} 

\begin{remark}
Our kind of \emph{the proliferation} of every $L_2$-orthogonal system is in the context of the Chumash, Bereishis, 26:12, \emph{Isaac sowed in the land, and in that year reaped a hunderdfold, thus had HASHEM blessed him}.
\end{remark} 

\section{Contributions to the theory of the Hardy-Littlewood integral that are generated by the Jacob's ladders} 

\subsection{} 

Let us remind that we have introduced Jacob's ladders 
\be \label{3.1} 
\vp_1(T)=\frac 12\vp(T) 
\ee 
in our paper \cite{2}, comp. also \cite{3}, where the function $\vp(T)$ is arbitrary solution of the nonlinear integral equation\footnote{Introduced also in \cite{2}.}
\be \label{3.2} 
\int_0^{\mu[x(T)]}\left|\zf\right|^2e^{\frac{-2}{x(T)}t}{\rm d}t, 
\ee 
where each admissible function 
\bdis 
\mu(y)\geq 7y\ln y
\edis  
generates a solution 
\be \label{3.3} 
y=\vp(T;\mu)=\vp(T). 
\ee 
We call the function $\vp_1(T)$ as Jacob's ladder since the analogy with the Jacob's dream in Chumash, Bereishis, 28:12. 

\subsection{} 

Next, 83 years after discovering the HLI formula (\ref{1.15}), we have proved new result, see \cite{2}. 

\begin{mydef91}
The Hardy-Littlewood integral \cite{1} 
\be \label{3.4} 
J(T)=\int_0^T\left|\zf\right|^2{\rm d}t
\ee 
has, in addition to previously known HLI formula (\ref{1.15}) possessing an unbounded error term at $T\to\infty$, the following infinite number of almost exact representations 
\be \label{3.5} 
\begin{split}
& \int_0^T\left|\zf\right|^2{\rm d}t=\vp_1(T)\ln\{\vp_1(T)\}+ \\ 
& (c-\ln 2\pi)\vp_1(T)+c_0+\mcal{O}\left(\frac{\ln T}{T}\right), 
\end{split}
\ee 
where $c_0$ is the constant from the Titschmarsh-Kober-Atkinson formula. 
\end{mydef91} 

\subsection{} 

Further, we have obtained also the following result\footnote{See \cite{2}, (6.2).}  

\begin{mydef92}
\be \label{3.6} 
T-\vp_1(T) \sim (1-c)\pi(T);\ \pi(T)\sim \frac{T}{\ln T},\ T\to\infty,  
\ee  
where the Jacob's ladder can be viewed as the complementary function to the function $(1-c)\pi(T)$ in the sense 
\be \label{3.7} 
\vp_1(T)+(1-c)\pi(T)\sim T,\ T\to\infty. 
\ee 
\end{mydef92} 

\begin{remark}
The following exchange 
\be \label{3.8} 
\frac{T}{\ln T}\to\pi(T),\ T\to\infty
\ee 
has been used in the first asymptotic formula (\ref{3.6}), see \cite{2}, (6.2). This, however, is completely correct in the asymptotic regions since the prime-number law. At the same time we can give also a kind of \emph{genetical motivation} for the exchange (\ref{3.8}). Namely 
\be \label{3.9} 
\zeta(s)=\exp\left\{s\int_2^\infty\frac{\pi(x)}{x(x^s-1)}{\rm d}x\right\},\ \re\{s\}>1, 
\ee 
together with the analytic continuation on $\mbb{C}\setminus \{1\}$. In this direction the prime-counting function $\pi(x)$ is one of five function which together generate the function $\zeta(s)$. 
\end{remark} 

\subsection{} 

Next, we have proved in \cite{7}, (3.4) the existence of almost linear increments of the Hardy-Littlewood integral (\ref{3.4}). 

\begin{mydef93}
For every fixed natural number $k$ and for every sufficiently big $T>0$ we have 
\be \label{3.10} 
\begin{split}
& \int_{\overset{r-1}{T}}^{\overset{r}{T}}\left|\zf\right|^2{\rm d}t=(1-c)\overset{r-1}{T}+\mcal{O}(T^{a+\delta}),\\ 
& 1-c\approx 0.42,\ r=1,\dots,k,\ T\to\infty, 
\end{split}
\ee 
where, see (\ref{1.17}), 
\be \label{3.11} 
\overset{r}{T}=\overset{r}{T}(T)=\vp_1^{-r}(T). 
\ee 
\end{mydef93} 

\begin{remark}
The existence of almost linear increments (linear in the main terms) in eq. (\ref{3.10}) is a new phenomenon in the theory of the Riemann's zeta-function on the critical line. It is true, for example\footnote{Comp. (\ref{1.16}).}, that  
\be \label{3.12} 
\int_{\overset{r-1}{T}}^{\overset{r}{T}}\left|\zf\right|^2{\rm d}t=(1-c)\overset{r-1}{T}+\mcal{O}(T^{\frac 13+\delta}),\ T\to\infty. 
\ee 
\end{remark} 

\begin{remark}
The formula (\ref{3.6}) represents the first interaction\footnote{See \cite{2}, (2.1), (2.5) and the section 6 ibid.} between the HLI formula\footnote{Comp. (\ref{1.15}).} and our formula, comp. (\ref{3.5}), while the formula (\ref{3.10}) is the result of the second interaction, see \cite{7}, (7.5) and (7.6), between the two mentioned formulas. 
\end{remark} 

\begin{remark}
It is clear that new contributions to the theory of the Hardy-Littlewood integral (1918) are represented also by: 
\begin{itemize}
	\item[(A)] the set of the thirteenth equivalents of the Fermat-Wiles theorem; 
	\item[(B)] the set of all corresponding functionals that generate the set in (A). 
\end{itemize}
\end{remark} 

\section{On almost exact decomposition of some increments of the Hardy-Littlewood integral by means of the Raabe's integral and Jacob's ladder} 

\subsection{} 

Let us remind the Raabe's integral 
\be \label{4.1} 
\int_a^{a+1}\ln\Gamma(t){\rm d}t=a\ln a - a + \ln\sqrt{2\pi},\ a>0 
\ee 
and also the formula 
\be \label{4.2} 
\begin{split}
& \int_0^{\overset{r}{T}}\left|\zf\right|^2{\rm d}t=\overset{r-1}{T}\ln\overset{r-1}{T}-(\ln 2\pi -c)\overset{r-1}{T}+c_0+ \\ 
& \mcal{O}\left(\frac{\ln T}{T}\right),\ t=1,\dots,k, 
\end{split}
\ee 
that holds true for every fixed natural number $k$ and which follows from almost exact formula (\ref{3.5}) by making the substitution\footnote{Comp. (\ref{1.18}).}
\bdis 
T\to\overset{r}{T};\ \vp_1(\overset{r}{T})=\overset{r-1}{T}. 
\edis 
Putting $a\to\overset{r-1}{T}$ in (\ref{4.1}) we immediately get obtain 
\be \label{4.3} 
\int_{\overset{r-1}{T}}^{\overset{r-1}{T}+1}\ln\Gamma(t){\rm d}t=\overset{r-1}{T}\ln \overset{r-1}{T}-\overset{r-1}{T}+\ln\sqrt{2\pi}, 
\ee  
and by (\ref{4.2}) and (\ref{4.3}) we get the result 
\be \label{4.4} 
\begin{split}
& \int_0^{\overset{r}{T}}\left|\zf\right|^2{\rm d}t-\int_{\overset{r-1}{T}}^{\overset{r-1}{T}+1}\ln\Gamma(t){\rm d}t= \\ 
& -(\ln 2\pi-c-1)\overset{r-1}{T}+c_0-\ln\sqrt{2\pi}+\mcal{O}\left(\frac{\ln T}{T}\right). 
\end{split}
\ee  
Next, the translation $r\to r+1$ in (\ref{4.4}) gives us 
\be \label{4.5} 
\begin{split}
& \int_0^{\overset{r+1}{T}}\left|\zf\right|^2{\rm d}t-\int_{\overset{r}{T}}^{\overset{r}{T}+1}\ln\Gamma(t){\rm d}t= \\ 
& -(\ln 2\pi-c-1)\overset{r}{T}+c_0-\ln\sqrt{2\pi}+\mcal{O}\left(\frac{\ln T}{T}\right). 
\end{split}
\ee 
And finally, by subtraction of the eq. (\ref{4.4}) from the eq. (\ref{4.5}) we have the following statement. 

\begin{mydef11}
It is true that 
\be \label{4.6} 
\begin{split}
& \int_{\overset{r}{T}}^{\overset{r+1}{T}}\left|\zf\right|^2{\rm d}t=\int_{\overset{r}{T}}^{\overset{r}{T}+1}\ln\Gamma(t){\rm d}t-
\int_{\overset{r-1}{T}}^{\overset{r-1}{T}+1}\ln\Gamma(t){\rm d}t- \\ 
& (\ln 2\pi-1-c)(\overset{r}{T}-\overset{r-1}{T})+\mcal{O}\left(\frac{\ln T}{T}\right),\ T\to\infty, 
\end{split}
\ee 
where 
\be \label{4.7} 
\ln 2\pi-c-1\approx 0.68,\ r=1,\dots,k
\ee 
for every fixed $k\in\mbb{N}$. 
\end{mydef11} 

Next,the summation of the formula (\ref{4.6}) gives the following result. 

\begin{mydef41}
\be \label{4.8} 
\begin{split}
	& \int_{\overset{1}{T}}^{\overset{k+1}{T}}\left|\zf\right|^2{\rm d}t=\int_{\overset{k}{T}}^{\overset{k}{T}+1}\ln\Gamma(t){\rm d}t-
	\int_{T}^{T+1}\ln\Gamma(t){\rm d}t- \\ 
	& (\ln 2\pi-1-c)(\overset{k}{T}-T)+\mcal{O}_k\left(\frac{\ln T}{T}\right),\ T\to\infty, 
\end{split}
\ee 
for every fixed $k\in\mbb{N}$. 
\end{mydef41} 

\begin{remark}
It is true by (\ref{4.6}): the areas $\mcal{A}_r^0$ generated, as usually, by the basic functions 
\be \label{4.9} 
y_r^0(t)=\left|\zf\right|^2,\ t\in[\overset{r-1}{T},\overset{r}{T}],\ r=1,\dots,k 
\ee 
are expressed by the big-oh exactness 
\be \label{4.10} 
\mcal{O}\left(\frac{\ln T}{T}\right),\ T\to\infty 
\ee 
by means of corresponding linear combination of areas $\mcal{A}_r^1,\mcal{A}_r^2,\mcal{A}_r^3$ that are generated by the following functions 
\be \label{4.11} 
\begin{split}
& y_r^1(t)=\ln\Gamma(t),\ t\in[\overset{r}{T},\overset{r}{T}+1], \\ 
& y_r^2(t)=\ln\Gamma(t),\ t\in[\overset{r-1}{T},\overset{r-1}{T}+1], \\ 
& y_r^3(t)=\ln 2\pi-1-c,\ t\in[\overset{r-1}{T},\overset{r}{T}],  
\end{split}
\ee 
respectively. In the last case, of course, 
\be \label{4.12} 
\overset{r}{T}-\overset{r-1}{T}=\int_{\overset{r-1}{T}}^{\overset{r}{T}}1{\rm d}t. 
\ee 
\end{remark} 

\section{New equivalent of the Fermat-Wiles theorem connected with the Raabe's integral and the Jacob's ladder} 

\subsection{} 

In this section we use the formulas that follow from the eq. (\ref{4.6}) by the substitution 
\be \label{5.1} 
\overset{r}{T}(T)=\rho_r,\ \overset{r-1}{T}(T)=\vp_1(\rho_r)
\ee 
for every fixed $r$\footnote{Comp. (\ref{1.18}).}. 

Let us notice the following about the substitution (\ref{5.1}): 
\begin{itemize}
	\item[(A)] Since $r$ is fixed one, we put in (\ref{5.1}) simply 
	\be \label{5.2} 
	\rho_r=\rho, 
	\ee 
	of course, 
	\be \label{5.3} 
	\{ T\to\infty \} \ \Leftrightarrow \ \{\rho\to\infty\}. 
	\ee 
	\item[(B)] Since\footnote{See (\ref{1.19}).} $\overset{r}{T}\sim T$, then   
	\be \label{5.4} 
	T\sim\rho 
	\ee 
	and, consequently, for the last member on the right hand side of (\ref{4.6}), we obtain 
	\be \label{5.5} 
	\mcal{O}\left(\frac{\ln T}{T}\right)=\mcal{O}\left(\frac{\ln \rho}{\rho}\right),\ \rho\to\infty. 
	\ee 
	\item[(C)] And finally also\footnote{See (\ref{1.19}) and (\ref{3.6}).} 
	\be \label{5.6} 
	\overset{r}{T}-\overset{r-1}{T}=\rho-\vp_1(\rho)\sim (1-c)\frac{\rho}{\ln\rho},\ \rho\to\infty. 
	\ee 
	This will be the basic formula used in present section. 
\end{itemize} 

\begin{mydef51}
\be \label{5.7} 
\begin{split}
& \int_\rho^{\overset{1}{\rho}}\left|\zf\right|^2{\rm d}t=\\ 
& \int_\rho^{\rho+1}\ln\Gamma(t){\rm d}t-\int_{\vp_1(\rho)}^{\vp_1(\rho)+1}\ln\Gamma(t){\rm d}t+\mcal{O}\left(\frac{\ln \rho}{\rho}\right),\ \rho\to\infty, 
\end{split}
\ee 
where 
\be \label{5.8} 
\overset{1}{\rho}=\vp_1^{-1}(\rho). 
\ee 
\end{mydef51} 

\subsection{} 

Next, by our almost linear formula\footnote{Comp. (\ref{3.10}).}
\be \label{5.9} 
\int_\rho^{\overset{1}{\rho}}\left|\zf\right|^2{\rm d}t=(1-c)\rho+\mcal{O}(\rho^{\frac 13+\delta})
\ee 
we get from (\ref{5.7}) the result 
\be \label{5.10} 
\begin{split}
& \int_\rho^{\rho+1}\ln\Gamma(t){\rm d}t-\int_{\vp_1(\rho)}^{\vp_1(\rho)+1}\ln\Gamma(t){\rm d}t= \\ 
& (1-c)\rho+\mcal{O}\left(\frac{\ln \rho}{\rho}\right). 
\end{split}
\ee 
Consequently, by making use of the substitution 
\be \label{5.11} 
\rho=\frac{x}{1-c}\tau,\ x>0;\ \{\rho\to\infty\} \ \Leftrightarrow \ \{\tau\to\infty\}
\ee 
in (\ref{5.10}) we get 
\be \label{5.12} 
\begin{split}
& \int_{\frac{x}{1-c}\tau}^{\frac{x}{1-c}\tau+1}\ln\Gamma(t){\rm d}t-\int_{\vp_1(\frac{x}{1-c}\tau)}^{\vp_1(\frac{x}{1-c}\tau)+1}\ln\Gamma(t){\rm d}t=\\
& x\tau+\mcal{O}_x\left(\frac{\tau}{\ln\tau}\right), 
\end{split}
\ee 
where the constant in the error term depends upon $x$. Now, it follows from the last eq. 

\begin{mydef52}
\be \label{5.13} 
\begin{split}
& \lim_{\tau\to\infty}\frac{1}{\tau}\left\{
\int_{\frac{x}{1-c}\tau}^{\frac{x}{1-c}\tau+1}\ln\Gamma(t){\rm d}t-\int_{\vp_1(\frac{x}{1-c}\tau)}^{\vp_1(\frac{x}{1-c}\tau)+1}\ln\Gamma(t){\rm d}t
\right\}=x,\ x>0.
\end{split}
\ee 
\end{mydef52} Next, by making use of the substitution 
\be \label{5.14} 
x\to\frac{x^n+y^n}{z^n},\ x,y,z,n\in\mbb{N},\ n\geq 3
\ee 
in (\ref{5.13}) we obtain the following statement. 

\begin{mydef53}
It is true that 
\be \label{5.15} 
\begin{split}
& \lim_{\tau\to\infty}\frac{1}{\tau}\left\{
\int_{\frac{x^n+y^n}{z^n}\frac{\tau}{1-c}}^{\frac{x^n+y^n}{z^n}\frac{\tau}{1-c}+1}\ln\Gamma(t){\rm d}t-\int_{\vp_1(\frac{x^n+y^n}{z^n}\frac{\tau}{1-c})}^{\vp_1(\frac{x^n+y^n}{z^n}\frac{\tau}{1-c})+1}\ln\Gamma(t){\rm d}t
\right\}= \\ 
& \frac{x^n+y^n}{z^n} 
\end{split}
\ee 
for every fixed Fermat's rational 
\bdis 
\frac{x^n+y^n}{z^n}. 
\edis 
\end{mydef53} 

Finally, the following result holds true. 

\begin{mydef12}
The condition 
\be \label{5.16} 
\begin{split}
& \lim_{\tau\to\infty}\frac{1}{\tau}\left\{
\int_{\frac{x^n+y^n}{z^n}\frac{\tau}{1-c}}^{\frac{x^n+y^n}{z^n}\frac{\tau}{1-c}+1}\ln\Gamma(t){\rm d}t-\int_{\vp_1(\frac{x^n+y^n}{z^n}\frac{\tau}{1-c})}^{\vp_1(\frac{x^n+y^n}{z^n}\frac{\tau}{1-c})+1}\ln\Gamma(t){\rm d}t
\right\}\not= 1
\end{split}
\ee 
on the class of all Fermat's rationals represent the thirteenth equivalent to the Fermat-Wiles theorem. 
\end{mydef12} 

\begin{remark}
Of course, it would be sufficient to localize the condition (\ref{5.16}) on the case 
\bdis 
1-\epsilon<\frac{x^n+y^n}{z^n}<1+\epsilon 
\edis 
for every positive, small and fixed $\epsilon$, comp. \cite{11}, sec. 7. 
\end{remark}

I would like to thank Michal Demetrian for his moral support of my study of Jacob's ladders.


\begin{thebibliography}{29}
\bibitem{1} 
G.H. Hardy, J.E. Littlewood, Contribution to the theory of the Riemann zeta-function and the theory of the distribution of Primes, Acta Math. 41 (1), 119 -- 196, (1918). 
\bibitem{2}
J. Moser,
`Jacob's ladders and almost exact asymptotic representation of the Hardy-Littlewood integral`,
Math. Notes 88, (2010), 414-422, arXiv: 0901.3937. 
\bibitem{3}
J. Moser,
`Jacob's ladders, the structure of the Hardy-Littlewood integral and some new class of nonlinear integral equations`,
Proc. Steklov Inst. 276 (2011), 208-221, arXiv: 1103.0359. 
\bibitem{4}
J. Moser, Jacob's ladders, reverse iterations and new infinite set of $L_2$-orthogonal systems generated by the Riemann $\zf$-function, arXiv: 1402.2098v1. 
\bibitem{5} 
J. Moser, Jacob's ladders, interactions between $\zeta$-oscillating systems and $\zeta$-analogue of an elementary trigonometric identity, Proc. Stek. Inst. 299, 189 -- 204, (2017). 
\bibitem{6} 
J. Moser, Jacob's ladders and vector operator producing new generations of $L_2$-orthogonal systems connected with the Riemann's $\zf$-function, arXiv: 2302.0750.v3 
\bibitem{7} 
J. Moser, Jacob's ladders, existence of almost linear increments of the Hardy-Littlewood integral and new types of multiplicative laws, arXiv: 2304.09267. 
\bibitem{8} 
J. Moser, Jacob's ladders, almost linear increments of the Hardy-Littlewood integral (1918) and their relations to the Selberg's formula (1946) and the Fermat-Wiles theorem, arXiv: 2306.07648v1. 
\bibitem{9} 
J. Moser, Jacob's ladders, almost linear increments of the Hardy-Littlewood integral (1918), the classical Dirichlet's sum of the divisors (1849) and their relationship to the Fermat-Wiles theorem, arXiv: 2312.12085. 
\bibitem{10} 
J. Moser, Jacob's ladders, almost linear increments of the Hardy-Littlewood integral (1918) and their relation to the Titchmarsh's sums (1934) and to the Fermat-Wiles theorem, arXiv: 2401.03781. 
\bibitem{11} 
J. Moser, Jacob's ladders, Hardy-Littlewood integral (1918) and new asymptotic functional equation for the Euler's $\Gamma$-function together with the tenth equivalent of the Fermat-Wiles theorem, arXiv: 2403.17522. 
\bibitem{12} 
J. Moser, Jacob's ladders, logarithmic modification of the Hardy-Littlewood integral (1918), Titchmarsh's $\Omega$-theorem (1928) and new point of contact with the Fermat-Wiles theorem, arXiv: 2406.02278. 
\end{thebibliography}
\end{document}